\documentstyle{amsart}

\newtheorem{thm}{Theorem}[section]
\newtheorem{prop}[thm]{Proposition}
\newtheorem{lem}[thm]{Lemma}
\newtheorem{cor}[thm]{Corollary}

\theoremstyle{definition}

\newtheorem{remark}[thm]{Remark}
\newtheorem{example}[thm]{Example}

\theoremstyle{remark}

\numberwithin{equation}{section}

\newcommand{\Id}{\operatorname{Id}}
\newcommand{\UD}{\operatorname{UD}}

\newcommand{\Span}{\operatorname{Span}}

\newcommand{\Spec}{\operatorname{Spec}}

\newcommand{\Hom}{\operatorname{Hom}}

\newcommand{\RMaps}{\operatorname{RMaps}}

\newcommand{\bbA}{{\Bbb A}}

\newcommand{\bbZ}{{\Bbb Z}}

\newcommand{\bbP}{{\Bbb P}}

\newcommand{\isomo}{\overset{\sim}{=}}
\newcommand{\isomoto}{\overset{\sim}{\to}}

 \newcommand{\brokrarr}{\vphantom{\to}\mathrel{\smash{{-}{\rightarrow}}}}

\newcommand{\lf}{\mathopen}

\let\r=\mathclose

\newcommand{\notdiv}{\mathrel{\not|}}

\newcommand{\GL}{{\operatorname{GL}}}

\newcommand{\PGL}{{\operatorname{PGL}}}
\newcommand{\PGLn}{{\operatorname{PGL}_n}}

\newcommand{\Char}{\operatorname{char}} 

\newcommand{\Gal}{\operatorname{Gal}}
\newcommand{\Galois}{\Gal}
\newcommand{\lra}{\longrightarrow}

\newcommand{\Mat}{{\operatorname{M}}}
\newcommand{\Mn}{\Mat_n}

\newcommand{\sqf}{\operatorname{sqf}}

\newcommand{\Sym}{{\operatorname{S}}}
\newcommand{\Or}{{\bold O}}

\newcommand{\tr}{\operatorname{tr}}

\let\to=\longrightarrow

\begin{document}
\title[Characteristic polynomials]
{Conditions satisfied by characteristic polynomials 
in fields and division algebras}
\author[Z. REICHSTEIN and B. YOUSSIN]
{Z. Reichstein and B. Youssin} 
\address{Department of Mathematics, Oregon State University,
Corvallis, OR 97331}
\thanks{Z. Reichstein was partially supported by NSF grant DMS-9801675}
\email{zinovy@@math.orst.edu}
\address{Department of Mathematics and Computer Science,
University of the Negev, Be'er Sheva', Israel\hfill\break
\hbox{{\rm\it\hskip\parindent Current mailing address\/}}: 
Hashofar 26/3, Ma'ale Adumim, Israel}
\email{youssin@@math.bgu.ac.il}

\begin{abstract}
Suppose $E/F$ is a field extension. We ask whether or not there exists
an element of $E$ whose characteristic polynomial has 
one or more zero coefficients in specified positions. 
We show that the answer is frequently ``no''. We also prove similar
results for division algebras and show that the universal 
division algebra of degree $n$ does not have an element of trace 0 and
norm 1.
\end{abstract}

\maketitle
\tableofcontents

\section{Introduction}
\label{sect1}

Let $E/F$ be a field extension of degree $n$ and $\det: E \lra F$ be the
norm function.
For $x \in E$, we define $\sigma^{(i)}(x)$ by
\begin{equation} \label{ci}
 \det(\lambda1_F -x) = \lambda^n + \sigma^{(1)}(x)\lambda^{n-1} + \cdots + 
\sigma^{(n-1)}(x) \lambda + \sigma^{(n)}(x) \; .
\end{equation}
In particular, $\sigma^{(1)}(x) = 
-\tr(x)$ and $\sigma^{(n)}(x) = (-1)^n \det(x)$. In the sequel, whenever
we write $\sigma^{(i)}(x)$, we shall always understand $i$ to be 
an integer between $1$ and $n$. If the reference to the extension $E/F$
is not clear from the context, we will sometimes write 
$\sigma_{E/F}^{(i)}(x)$ in place of   $\sigma^{(i)}(x) \in F$.

If $A$ is a central simple algebra of degree $n$ with center $F$ then
we can define $\sigma^{(i)} =
\sigma_{A/F}^{(i)}$ in the same way. Here $\det$ 
in formula~\eqref{ci} should be intepreted as the reduced norm in 
$A \otimes_F F(\lambda)$. 

A number of interesting results, both in the theory of polynomials and
in the theory of central simple algebras, can be stated in terms of 
the existence (or nonexistence) of nontrivial solutions to 
systems of equations of the form
\begin{equation} \label{e.system}
 \sigma^{(i)}(x) = 0 \quad \text{\rm for} \quad i = i_1, \dots, i_r  \; .
\end{equation}
\begin{example} \label{ex1.1}
(Hermite~\cite{hermite}, Joubert~\cite{joubert}; see also Coray~\cite{coray})
If $E/F$ is a field extension of degree $5$ or $6$ and 
$\Char(F) \neq 3$ then there exists an element $x \in E$ such 
that $E=F(x)$ and $\sigma^{(1)}(x) = \sigma^{(3)}(x) = 0$. 
 
In classical language, this means that for $n = 5$ or $6$ every polynomial 
$f(t) = t^n + a_1t^{n-1} + \dots + a_n \in F[t]$ can be reduced,
via the Tschirnhaus transformation $t \mapsto x$, to the form 
$f(t) = t^n + b_1t^{n-1} + \dots + b_n \in F[t]$ with $b_1 = b_3 = 0$;
for details we refer the reader to~\cite{br}. 
\end{example}

\begin{example} \label{ex1.2}
Let $A$ be a central simple algebra of degree $n$ whose center
contains a primitive $n$th root of unity. Then $A$ is cyclic iff
there exists an element $x$ such that 
\[ \sigma^{(1)}(x) = \dots = \sigma^{(n-1)}(x) = 0 \; . \]
A conjecture of Albert asserts that every $A$ of prime (or, equivalently,
square-free) degree is cyclic. This conjecture is known to be 
true for $n = 2$, $3$ and $6$ (see~\cite[Section 3.2]{rowen-pi});
the remaining cases are open. 
\end{example}

\begin{example} \label{ex1.3}
(Haile~\cite{haile}; see also Brauer~\cite[Proposition~7.1.43]{rowen-rt}) 
Suppose $A$ is a central simple algebra of degree $n$ with center $F$.
Then there exists an $(n-1)$-dimensional $F$-subspace $W$ of $A$ such that
$\sigma^{(1)}(x) = \sigma^{(n-1)}(x) = 0$ for any $x \in W$.
\end{example}

\begin{example} \label{ex1.4}
(Rowen~\cite[Corollary 5]{rowen-brauer}) 
If $A$ is a central simple algebra of odd degree with center $F$ then
there exists an element $x \in A - \{ 0 \}$ such that $\sigma^{(1)}(x) =
\sigma^{(2)}(x) = 0$. 

Note that if $\Char(F) \neq 2$, this follows easily from a theorem 
of Springer (see e.g.,~\cite[Remark 14.3]{reichstein}); however,
the above result is true even if $\Char(F) = 2$.
\end{example}

In~\cite{reichstein} the first author showed that in many cases
equations of the form $\sigma^{(i)}(x) = 0$ or $\tr(x^i) = 0$ 
and systems of the form 
$\sigma^{(1)}(x) = \sigma^{(i)}(x) = 0$ or $\tr(x) = \tr(x^i) = 0$
do not have nontrivial solutions. In particular,
the theorem of Hermite and Joubert, cited in Example~\ref{ex1.1},
fails for field extensions of degree $n = 3^m$ or $3^m + 3^l$, 
with $m > l \geq 0$.  In this paper we revisit this subject from a more
geometric point of view. 

\subsection*{Notational conventions}
Throughout this paper $n$ will denote the degree of the field extension 
or division algebra we are considering, and $\sqf(n)$ will denote 
the square-free part of $n$. We will always work over a fixed ground field $k$.

Let $K$ be a field containing a primitive $r$th root of unity
$\zeta_r$ (in particular, we assume that $r$ is prime to $\Char(K)$),
and let $z, w \in K$.  Recall that a symbol algebra
$(z, w)_r$ is defined as 
\begin{equation} \label{e.symbol}
(z, w)_r = K\{x, y\}/( x^r = z \, , \;  y^r = w \, , \;  yx = \zeta_r xy ) 
\, ;
\end{equation}
cf.~\cite[p. 194]{rowen-rt}.  We now define the algebra $D_n$ as follows.
Write $n = p_1 \dots p_s$ as a product of (not necessarily distinct) primes.
Let $K = k(z_1, w_1, \dots , z_s, w_s)$, where 
$z_1, w_1, \dots, z_s, w_s$ are independent variables over $k$ and let
\begin{equation} \label{e.D} 
D_n = (z_1, w_1)_{p_1} \otimes_K \dots \otimes_K (z_s, w_s)_{p_s} \; .
\end{equation}
Note that $D_n$ is a division algebra of degree $n$ and 
exponent $\sqf(n)$, with center $K$. 

Finally recall that the {\em universal division 
algebra} $\UD(n)$ is the subalgebra of $\Mn(k(s_{ij}, t_{ij}))$ generated,
as a division algebra, by two generic $n \times n$-matrices $(s_{ij})$
and $(t_{ij})$. Here $s_{ij}$ and $t_{ij}$ are $2n^2$ independent 
variables over $k$. For details of this construction, see, 
e.g.,~\cite[Section 3.2]{rowen-pi}.

\subsection*{Main results}

\begin{thm} \label{thm3} 
Suppose $\Char(k) \notdiv n!$ and $D = D_n$ or $\UD(n)$. Then the system 
 \begin{equation} \label{e.thm3}
 \left\{ \begin{array}{l} \sigma^{(i)}(x_1) = \dots = \sigma^{(i)}(x_m) \\
       \sigma^{(j)}(x_1 \dots x_m) = 0 \end{array} \right. 
 \end{equation}
has no nontrivial solutions in $D$, provided that
$i$ and $m$ are divisible by $\sqf(n)$. 
\end{thm}

Here, as usual, a solution $(x_1, \dots, x_s)$ is trivial 
if $x_1 = \dots = x_s = 0$ and nontrivial otherwise. Note that
the assertion of the theorem for $\UD(n)$ is a formal consequence of
the assertion for $D_n$, because of the specialization property 
of $\UD(n)$.  However, our proof will treat the two cases in parallel,
since both are proved by the same argument. Theorem~\ref{thm3} 
can be generalized in several directions; some generalizations
are discussed at the end of Section~\ref{sect3}. 

We now record three consequences of Theorem~\ref{thm3}, which we feel
deserve a special mention.

\begin{cor} \label{cor4} Suppose $\Char(k) \notdiv n!$, $D = D_n$ 
or $\UD(n)$, and $m$ is divisible by $\sqf(n)$. 

(a) $\sigma^{(m)}(x) \neq 0$ for any $x \in D - \{ 0 \}$.

\smallskip
(b) If $\det(x_1) = \dots = \det(x_m)$ for some $x_1, \dots, x_m \in
D - \{ 0 \}$ then $\tr(x_1 \dots x_m) \neq 0$. 

(c) $D$ does not have an element 
of (reduced) norm 1 and (reduced) trace 0. 
\end{cor}

To prove part (a), we assume the contrary and substitute $i = m$, $x_1 = x$ 
and $x_2 = \dots = x_m = 0$ into~\eqref{e.thm3} to obtain a contradiction. 
To prove part (b), we apply Theorem~\ref{thm3} with $i = n$ and $j = 1$.
Finally, if $\det(x) = 1$ then setting $x_1 = x$ and $x_2 = \dots = x_m = 1$ 
in part (b), we obtain $\tr(x) \neq 0$, thus proving part (c).
\qed

\smallskip
The commutative counterpart of the universal division algebra is the
{\em general field extension} $L_n/K_n$ defined as follows:
\begin{eqnarray} \label{general} K_n = k(a_1, \dots, a_n) \; \, 
\text{\rm and} \; \, 
L_n = K_n[x]/(x^n + a_1 x^{n-1} + \dots + a_n) \, , 
\end{eqnarray}
where $a_1, \dots, a_n$ are algebraically independent indeterminates 
over $k$. 

\begin{thm} \label{thm1} Let $n_1$ and $n_2$ be positive integers,
and $L_n/K_n$ be the general field extension of degree $n= n_1 + n_2$.
Then the system of equations 
\begin{equation} \label{e.thm1} 
\tr(x^{m_1}) = \tr(x^{m_2}) = 0 
\end{equation}
has no nontrivial solutions $x \in L_n^*$, provided that 

\smallskip
(i) $n_1 n_2 \neq 0$ and $(-\frac{n_2}{n_1})^{m_2-m_1} \neq 1$ in $k$.
%
%

\smallskip
(ii) each $\sqf(n_i)$ ($i = 1, 2$) divides $m_1$ or $m_2$ (and possibly both).
\end{thm}

Note that if $\Char(k) =0$ then condition (i) holds unless $m_1 = m_2$
or $n_1 = n_2$ and $m_2 - m_1$ is even. If we replace (i) by 
a more complicated condition, we can also show that the system
$\sigma^{(m_1)}(x) = \sigma^{(m_2)}(x) = 0$ has no nontrivial solutions;
see Section~\ref{sect5}.

It is interesting to note that Theorem~\ref{thm3} 
and Corollary~\ref{cor4} remain true if $D$ is replaced by $L_n$; 
see Remark~\ref{rem.thm3-fld}. On the other hand, Theorem~\ref{thm1} 
fails if $L_n$ is replaced by $\UD(n)$; see Remark~\ref{rem.thm1-div}. 

All of the main results in this paper are proved by the same general method.
based on the {\em Going Down Theorem}~\ref{going-down}. This method is
outlined in Section~\ref{sect2.1}. In particular, our proofs 
of Theorems~\ref{thm3} and~\ref{thm1}, given in
Sections~\ref{sect3} and~\ref{sec:pf.thm1}, are applications of
Propositions~\ref{prop:gen.div-alg} and~\ref{prop:down.Sn}
respectively.
Proposition~\ref{prop:gen.div-alg} says a system 
of equations, such as~\eqref{e.thm3}, has no nontrivial 
solutions in a ``sufficiently generic'' division algebra 
if a certain projective $\PGLn$-variety, constructed from 
this system, does not have $H$-fixed points for some abelian 
subgroup $H$ of $\PGLn$. Proposition~\ref{prop:down.Sn}
gives a similar criterion for nonexistence of solutions in field 
extensions. Other applications of this approach and some generalizations 
are presented in Sections~\ref{sect5}--\ref{sect7}.



\subsection*{Acknowledgements} The authors would like to thank 
A. R. Wadsworth for helpful discussions.

\section{The Going Down Theorem and its applications}
\label{sect2.1}

The following result will play a key role in
the sequel. A simple proof, due to Koll{\'a}r and Szab{\'o},
can be found in~\cite[Appendix]{ry1}.  Assume that $k$ is an
algebraically closed base field, and that all varieties, group actions 
and maps are defined over $k$.

\begin{thm}[The Going Down Theorem] \label{going-down} 
Let $H$ be a finite abelian group acting
on algebraic varieties $X$ and $Y$ and let $f \colon X \brokrarr Y$ 
be an $H$-equivariant rational map.  If $X$ has a smooth $H$-fixed point 
and $Y$ is projective then $Y$ has an $H$-fixed point.
\qed
\end{thm}

\subsection*{$\protect\Sym_n$-varieties}

Let $L/K$ be a separable field extension of degree $n$, let $L'$ be
the normal closure of $L$ over $K$, and 
$\Gal(L'/K)=G$. Note that $G$ acts on the set of embeddings
$L \hookrightarrow L'$ and thus is naturally realized as
a transitive subgroup of $\Sym_n$. For each $i = 1, \dots, n$
choose $g_i \in \Sym_n$ such that $g_i(1) = i$.
The embedding of $G$ in $\Sym_n$
defines a (permutation) action of $G$ on $\bbA^n$
and thus a diagonal actions on $(\bbA^n)^m$ for every $m \geq 1$.

Let $P(x_{11}, \dots, x_{1n}; \dots; x_{m1}, \dots, x_{mn})
 \in k[(\bbA^n)^m]$ be a $G$-invariant polynomial and let $a_1, \dots, a_m
\in L$.  
Then we can define $P(a_1,\dots,a_m)$ as $P(a_{11}, \dots, a_{1n}; \dots;
a_{m1}, \dots, a_{mn})$, where $a_{ij} = g_j(a_i) \in L'$. A priori,
$P(a_1, \dots, a_m) \in L'$; however, since $P$ is
$G$-invariant polynomial, $P(a_1, \dots, a_m)$ actually lies in
$(L')^G = K$.

In the sequel we shall assume that
$K$ is finitely generated over $k$ (and hence, so are $L$ and $L'$).

\begin{prop} \label{prop:down.Sn}
Let $Y$ be the subvariety of $\bbP((\bbA^n)^m)$ given by
$G$-invariant homogeneous polynomial equations $P_1=\dots=P_s=0$.
Suppose that $Y$ does not have $H$-fixed points for some abelian
subgroup $H\subset G$.
Assume that there exists a $G$-variety $X$ which has a smooth
$H$-fixed point and such that $k(X)=L'$ as fields with $G$-action.
Then the system of equations
\begin{equation} \label{eqn:gen.in.field}
P_1(a_1,\dots,a_m)=\dots=P_s(a_1,\dots,a_m)=0
\end{equation}
has no nontrivial solutions in $L$.
\end{prop}

We remark that if $\Char(k)=0$ then a $G$-variety
$X$ such that $k(X)=L'$ (as $G$-fields) 
always exists; see~\cite[Proposition 8.6 and Example 8.4c]{tg}. 
Moreover, we can choose $X$ to be smooth and projective;
see~\cite[Proposition 2.2]{ry2}. In view of Theorem~\ref{going-down}, 
the presence of an $H$-fixed point on such an $X$ is a birational 
invariant, i.e., is independent of the choice of the (smooth projective)
model. 

\begin{pf}
Suppose $(a_1,\dots,a_m)\in L^m\subset k(X)^m$ is a non-trivial solution 
of \eqref{eqn:gen.in.field} and let $a_{i1}, \dots, a_{in}$ be the
conjugates of $a_i$ in $L'$. Then 
\[ f \colon x \mapsto [a_{11}(x): a_{12}(x): \dots : a_{mn}(x)] \]
is a $G$-equivariant rational map $X\brokrarr\bbP((\bbA^n)^m)$.
By our choice of $a_1,\dots,a_m$, the image of $f$ lies in $Y$.
Applying Theorem~\ref{going-down} to the rational map
$f \colon X \brokrarr Y$, we conclude that $Y$ has an $H$-fixed point, 
a contradiction.
\end{pf}

In the sequel we shall use use Proposition~\ref{prop:down.Sn}
only for $m=1$; the statement for general $m$ is intended
to make it parallel to Proposition~\ref{prop:gen.div-alg} below.

\subsection*{$\protect\PGLn$-varieties} 

Let $P\in k[(\Mn)^m]^\PGLn$; it is a polynomial in the entries of
$m$ matrices $U_1,\dots,U_m$ invariant under simultaneous conjugation.
If $A$ is a central simple algebra of degree $n$ and $a_1, \dots, a_m
\in A$ then we can define $P(a_1, \dots, a_m)$ as follows.
Split $A$ by the algebraic closure $\overline{K}$ of $K$:
$A\otimes_K\overline K\simeq\Mn(\overline K)$.
Thus $A\hookrightarrow \Mn(\overline K)$, and we can evaluate
$P(a_1,\dots,a_m)\in\overline K$.

\begin{lem} \label{lem:inv.polyn}
$P(a_1,\dots,a_m)$ lies in $K$ and is independent of the choice of
the isomorphism $A\otimes_K\overline K\simeq\Mn(\overline K)$.
\end{lem}

\begin{pf}
Any two choices of the isomorphism 
$A\otimes_K\overline K\simeq\Mn(\overline K)$ differ by  
conjugation by some $g \in \PGLn(\overline K)$. 
Since $P$ is $\PGLn$-invariant, conjugation by $g$
does not change the value of $P(a_1,...,a_m)$. 

Consider the action of $\Gal(\overline K/K)$ on $\Mn(\overline K)$;
for any $\sigma\in\Gal(\overline K/K)$ and
$B_1,\dots,B_n\in\Mn(\overline K)$,
$P(\sigma(B_1),\dots,\sigma(B_n))=\sigma(P(B_1,\dots,B_m))$.
The composition
\[
\Mn(\overline K)\isomoto A\otimes_K\overline K@>\Id\otimes\sigma^{-1}>>
A\otimes_K\overline K\isomoto\Mn(\overline K)
@>\sigma>>\Mn(\overline K)
\]
is an automorphism of $\Mn(\overline{K})$ whose restriction to the center
$\overline K$ is trivial.  Hence, this composition is given by
conjugation by some $g \in \PGLn(\overline K)$.
It follows that for $a_1,\dots,a_m\in A$,
$P(a_1,\dots,a_m)$ is fixed by $\Gal(\overline K/K)$ and
thus lies in $K$.
\end{pf}

Note that the Lemma is an immediate consequence of the fact that
$k[(\Mn)^m]^{\PGLn}$ is generated by elements of the form
$\sigma^{(i)}(U)$, where $U$ is a monomial in the $m$-matrices 
$U_1, \dots, U_m$. The latter was proved by Sibirskii~\cite{sibirskii}
and Procesi~\cite{procesi1} in the case $\Char(k) = 0$ and, more recently,
by Donkin~\cite{donkin} in prime characteristic. The elementary argument
given above allows us to avoid appealing to this more difficult result.

Next we recall that is $F$ be a finitely generated field extension 
of $k$ then an element of $H^1(F, \PGLn)$ may be interpreted either as
a central simple algebra $D$ of degree $n$ with center $F$ or, alternatively,
as a generically free $\PGL_n$-variety $X$ such that $k(X)^{\PGLn} = F$.
It is shown in~\cite{reichstein} (under the assumption $\Char(k) = 0$) that
$D \isomo \RMaps_{\PGLn}(X, \Mn)$ = the algebra of $\PGLn$-equivariant rational
maps from $X$ to $\Mn$; see also~\cite[Section 3]{ry2}.
Note that the above isomorphism is an isomorphism of $F$-algebras,
where we identify $f \in F = k(X)^{\PGLn}$
with the $\PGLn$-equivariant rational map $X \brokrarr \Mn(k)$
given by $x \mapsto f(x)I_n$. (Here $I_n$ denotes the $n \times n$-identity 
matrix.)

\begin{prop} \label{prop:gen.div-alg}
Let $Y$ be the subvariety of $\bbP((\Mn)^m)$ cut out by $\PGLn$-invariant 
homogeneous polynomial equations $P_1=\dots=P_s=0$. Suppose
$Y$ has no fixed points for some
finite abelian subgroup $H$ of $\PGLn$.
Then the system of equations
\begin{equation} \label{eqn:gen.in.Mn^m}
P_1(x_1,\dots,x_m)=\dots=P_s(x_1,\dots,x_m)=0
\end{equation}
has no nontrivial solutions in any central simple algebra $D$ of the form
$D=\RMaps_{\PGLn}(X,\Mn)$, where $X$ is a generically free 
$\PGLn$-variety which has a smooth $H$-fixed point.
\end{prop}

\begin{pf}
Suppose the system \eqref{eqn:gen.in.Mn^m} has a nontrivial solution
$(x_1,\dots,x_m)$.
As $D=\RMaps_\PGLn(X,\Mn)$, each $x_i$ can be interpreted as a rational
$\PGLn$-invariant map $X\brokrarr\Mn$; collectively, these elements
define a rational
$\PGLn$-equivariant map $f \colon X \brokrarr \bbP((\Mn)^m)$.
By our choice of $x_1, \dots, x_m$, the image of this map lies in $Y$.
By Theorem~\ref{going-down}, $Y$ has a $H$-fixed point, a
contradiction.
\end{pf}

\section{Abelian subgroups}
\label{sec:ab.sgp}

In order to use Propositions~\ref{prop:down.Sn} 
and~\ref{prop:gen.div-alg}, we need a description of abelian 
subgroups $H$ of $\Sym_n$ and $\PGLn$.  In this section we introduce
the abelian subgroups that will be used in subsequent applications.

We shall assume that the base field $k$ contains all roots of unity.
For a finite abelian group $A$ 
of order prime to $\Char(k)$, we shall
denote its dual group $\Hom(A,k^*)$ by $A^*$.

\subsection*{Abelian subgroups of $\protect\Sym_n$}
Let $A = \{ a_1, \dots, a_n \}$ be an abelian group of order $n$. 
The right multiplication action of
$A$ on itself gives rise to an embedding 
\[ \psi_A \colon A \hookrightarrow \Sym_n \; . \]
(Note that if we relabel the elements of $A$, $\psi_A$ will change by
an inner automorphism of $\Sym_n$.)
Given a character $\chi \colon A \lra k^*$, let 
\[ R_{\chi} = (\chi(a_1), \dots, \chi(a_n)) \; . \]
It is easy to see that $ k^n = \bigoplus_{\chi \in A^*} 
\Span_k(R_{\chi})$ 
is a decomposition of $k^n$ as a direct sum of 1-dimensional
character spaces for the permutation action of $A$ on $k^n$ (via $\psi_A$);
moreover, the character associated to $\Span_k(R_{\chi})$ is precisely 
$\chi^{-1}$.

In the sequel we will be interested in the permutation action of 
\begin{equation} \label{e.product-group}
 H = \psi_{A_1}(A_1) \times \psi_{A_2}(A_2) \subset \Sym_{n_1} \times
\Sym_{n_2} \subset \Sym_n 
\end{equation}
on $k^n$. Here $A_1$ and $A_2$ are abelian groups of order $n_1$ and $n_2$ 
respectively and $n = n_1 + n_2$. For future reference, 
we decompose this action as a direct sum of character spaces. We shall write
elements of $k^n = k^{n_1 + n_2}$ as $(R', R'')$, where $R' \in k^{n_1}$ and
$R'' \in k^{n_2}$. Let $V_0 =
\{ (\underbrace{a,\dots,a}_{\text{$n_1$ times}},
\underbrace{b,\dots,b}_{\text{$n_2$ times}}) \mid
a, b \in k \}$.

\begin{lem} \label{lem.subgr-S_n}
\[\textstyle
k^n = 
V_0\oplus\left( \bigoplus_{\chi \in A_1^*} \Span_k(R_{\chi}, 0)\right)
\oplus\Bigl(\bigoplus_{\eta \in A_2^*} \Span_k(0, R_{\eta})\Bigr)
\]
is a decomposition of 
$k^n$ as a direct sum of character spaces for the $H$-action defined above.
Here $V_0$ is a 2-dimensional subspace with trivial associated
character; the remaining $n-2$ summands are 1-dimensional 
subspaces with distinct nontrivial characters.
\end{lem}

\begin{pf}
The proof of this lemma amounts to verifying that the summands
of the above decomposition are, indeed, character spaces and finding
their characters. We leave the details of the reader.
\end{pf}

\subsection*{Abelian subgroups of $\protect\PGLn$}

Let $A$ be an abelian subgroup of order $n$ and $V = k[A]$.
The group $A$ acts on $V$ by the regular 
representation $a \mapsto P_a \in \GL(V)$, where
\[\textstyle
 P_a(\sum_{b \in A} c_b b) = \sum_{b \in A} c_b ab  \]
for any $a \in A$ and $c_b \in k$. 
The dual group $A^{\ast}$ acts on $V$ by the representation
$\chi \mapsto D_{\chi} \in \GL(V)$, where
\[\textstyle
 D_{\chi}( \sum_{a \in A} c_a a) = \sum_{a \in A} c_a \chi(a) a  \]
for any $\chi \in A^{\ast}$ and $c_a \in k$. 
Note that in the basis $\{ a  \, | \, a \in A \}$ of $V$,
each $P_{a}$ is represented by a permutation matrix and each $D_{\chi}$ is
represented by a diagonal matrix; this explains our choice of the letters
$P$ and $D$.  It is easy to see that
\begin{equation} \label{e.pgln}
D_{\chi} P_a = \chi(a) P_a D_{\chi} \; ;
\end{equation}
hence, we have constructed an embedding
\begin{equation} \label{e.phi_A}
 \phi_A \colon A \times A^* \hookrightarrow \PGL(V) = \PGLn 
\end{equation}
given by $(a,  \chi) \mapsto \overline{P_a} \cdot \overline{D_{\chi}}$, 
where $\overline{P_a}$ and $\overline{D_{\chi}}$ are the elements 
of $\PGL(V)$, represented, respectively, by $P_a$ and $D_{\chi} \in \GL(V)$.

For future reference we record two simple lemmas.

\begin{lem} \label{lem.subgr-PGL_n}
For each $a \in A$ and $\chi \in A^*$, $V_{a, \chi} = \Span_k(P_a D_{\chi})$
is a 1-dimensional $H$-invariant subspace of $\Mn$, with associated character
$(b, \eta) \mapsto \chi^{-1}(b) \eta(a)$. Moreover, the $n^2$ matrices
$P_aD_{\chi}$ form a $k$-basis of $\Mn$. 
\end{lem}

\begin{pf} The first assertion is immediate from~\eqref{e.pgln}. 
Since the $n^2$ characters associated to the spaces $V_{a, \chi}$ are
distinct, the second assertion now follows from linear independence 
of characters. 
\end{pf}

\begin{lem} \label{lem.P_aD_chi}
Let $A$ be an abelian group of order $n$ and $(a, \chi)$ be an element 
of order $c$ in $A \times A^*$. 

\smallskip 
(a) $(P_aD_{\chi})^c = \epsilon I_n$, where 
$\epsilon = \chi(a)^{\frac{1}{2}c(c-1)} = \pm 1$
and $I_n$ is the $n \times n$-identity matrix.

\smallskip 
(b) The characteristic polynomial of $P_aD_{\chi}$ is $r(t) = 
(t^c -\epsilon)^{\frac{n}{c}}$. 

\smallskip
(c) Assume $\Char(k) \notdiv (\frac{n}{c})!$.
Then $\sigma^{(i)}(P_a D_{\chi}) \neq 0$ for any $i$ divisible by $c$.
\end{lem}

\begin{pf} (a) The identity
$(P_aD_{\chi})^c = \epsilon I_n$, where $\epsilon =
\chi(a)^{\frac{1}{2}c(c-1)}$, is immediate from~\eqref{e.pgln}. 
To see that $\epsilon = 1$ or $-1$, note that $\epsilon^2 =
(\chi(a)^c)^{c-1} = 1^{c-1} = 1$. 

\smallskip
(b) Let $C$ be the cyclic subgroup of $A \times A^*$
generated by $(a, \chi)$, so that $c = |C|$.
For each $\alpha \in (A \times A^*)/C$, let
$V_{\alpha}$ be the vector subspace of $\Mn$ spanned by 
 $(b, \eta) \in \alpha$. Each $V_{\alpha}$ is a $c$-dimensional
 subspace of $\Mn$, which is stable under right multiplication 
by $P_aD_{\chi}$. Since the matrices $P_aD_{\chi}$ form a basis
of $\Mn$ as $(a, \chi)$ ranges over $A \times A^*$ 
(see Lemma~\ref{lem.subgr-PGL_n}), we can write
\begin{equation} \label{e.V_alpha} \textstyle
\Mn = \bigoplus_{\alpha \in (A \times A^*)/C} V_{\alpha} \; .
\end{equation}
By part (a), $(P_aD_{\chi})^c = \epsilon I_n$. It is now easy to see
that the characteristic polynomial for the action 
of $P_aD_{\chi}$ on each $V_{\alpha}$ is $p(t) = t^c - \epsilon$. 
Consequenly, the charactersistic polynomial for the left multiplication 
action of  
$P_aD_{\chi}$ on $\Mn$ is $q(t) = p(t)^{n^2/c}$ (one factor of $p(t)$ for
each subspace $V_{\alpha}$ in~\eqref{e.V_alpha}), and the characteristic
polynomial of the $n \times n$-matrix $P_aD_{\chi}$ (or, equivalently, of
its action on $n \times 1$-column vectors) is 
\[ r(t) = q(t)^{1/n} = p(t)^{n/c} = (t^c - \epsilon)^{n/c} \ , \] 
as claimed.

\smallskip
(c) The binomial formula 
tells us that under our assumption on $\Char(k)$,
every monomial of the form $t^{n-i}$ with $i$ divisible by $c$ 
(and $i \leq n$), appears in $r(t)$ with a nonzero coefficient. 
In other words, for these values of $i$, 
$\sigma^{(i)}(P_a D_{\chi}) \neq 0$, as claimed.
\end{pf}

\section{Proof of Theorem~\ref{thm3}}
\label{sect3}
 
We may (and will, throughout this section)
assume without loss of generality that $k$ is
an algebraically closed field. Otherwise we can simply replace $D$ by
$\overline{D} = D \otimes_k \overline{k}$, 
where $\overline{k}$ is the algebraic closure of $k$: if 
the system~\eqref{e.thm3} has no nontrivial solutions in $\overline{D}$, 
it cannot have one in $D$.

Our goal is to deduce Theorem~\ref{thm3} as a special case of
Proposition~\ref{prop:gen.div-alg}. We shall now proceed
to introduce the finite abelian group $H$ and the 
$\PGLn$-varieties $X$ and $Y$ and to show that they satisfy 
the conditions of Proposition~\ref{prop:gen.div-alg}. We will then apply
Proposition~\ref{prop:gen.div-alg} with these $H$, $X$, and $Y$,
to conclude that the system~\eqref{e.thm3} has no nontrivial solutions 
in $D_n$ or $\UD(n)$. 

\subsection*{The group $H$} We define $H$ to be the finite
abelian subgroup of $\PGLn$ given by 
\begin{eqnarray} \label{e.H}
H = A \times A^* \stackrel{\phi_{A}}{\hookrightarrow} \PGLn \, ,
& \text{\rm where} &
A = \bbZ/p_1\bbZ \times \dots \times \bbZ/p_s\bbZ \, . 
\end{eqnarray}
Here, as in Section~\ref{sect1},
$n = p_1 \dots p_s$, where $p_1, \dots, p_s$ are 
not necessarily distinct primes;
the inclusion $\phi_A$ is as in~\eqref{e.phi_A}.
Note that the assumption $\Char(k)\notdiv n!$ of Theorem~\ref{thm3}
implies that $|H|=n^2$ is prime to $\Char(k)$.

\subsection*{The variety $X$} 

We shall now write the algebras that come up in the statement 
of Theorem~\ref{thm3}, namely $D = \UD(n)$ and $D = D_n$, in the form
$\RMaps_{\PGLn}(X, \Mn)$ for specific $\PGLn$-varieties $X$. Note that
we do not assume $\Char(k) = 0$.

\begin{lem} \label{lem.X_D-a}
(Procesi) $\UD(n) = \RMaps_{\PGLn}(X, \Mn)$, where $X = (\Mn)^2$ and
$\PGLn$ acts on $X$ by simultaneous conjugation.
\end{lem}

\begin{pf} See~\cite[Theorem~14.16]{saltman}, cf.\
also~\cite[Theorem~2.1]{procesi2} or~\cite[Example~3.1]{ry2}.
\end{pf}

Let $G$ be an algebraic group, $S$ be a closed subgroup of $G$, and $Y$
be an affine $S$-variety. The groups $S$ and $G$ act on  $G \times Y$ via
respectively, $s(g, y) = (gs^{-1}, sy)$ and $g'(g, y) = (g'g, y)$; moreover, 
the two actions commute. Thus the quotient $(G \times Y)//S =
\Spec(k[G \times Y]^S)$ is a $G$-variety; we will denote it by
$G \ast_S Y$. We will restrict our attention to the case where $S$ is 
a finite group of order prime to $\Char(k)$. In this case a theorem
of Hilbert and Noether (see, e.g.,~\cite[Theorem~1.1]{smith})
tells us that $k[G \times Y]^S$ is a finitely generated 
$k$-algebra, i.e., $G*_S Y$ is again an affine variety (of finite type).

\begin{lem} \label{lem.X_D-b}
There exists a faithful $2s$-dimensional
linear representation $V$ of $H$ such that
$D_n \simeq \RMaps_{\PGLn}(X, \Mn)$, where $X = \PGLn \ast_{H} V$.
\end{lem}

\begin{pf} Choose a set of generators $a_1, \dots, a_s$ for $A$ and 
a ``dual'' set of generators $\chi_1, \dots, \chi_s$ for $A^*$ so that 
\[ \chi_i(a_j) = \left\{ \begin{array}{ll} 1 & \text{\rm if $i \neq j$} \\
                         \zeta_{p_i} & \text{\rm if $i = j$} \, , 
\end{array} \right. \]
where $\zeta_{p_i}$ is the same primitive $p_i$th root of unity 
used in defining $(z_i, w_i)_{p_i}$; see~\eqref{e.symbol}
and~\eqref{e.D}.
Consider the
faithful action of $H = A \times A^*$ on $V = k^{2s}$ given by
\begin{eqnarray*}
\lefteqn{(a, \chi) : (\alpha_1, \dots, \alpha_s, 
\beta_1, \dots, \beta_s)  \mapsto}  \nonumber \\ 
 & & (\chi^{-1}(a_1) \alpha_1, \dots, \chi^{-1}(a_s) \alpha_s, 
\chi_1(a)\beta_1, \dots, \chi_s(a)\beta_s) \, .
\end{eqnarray*}
Set $X = \PGLn *_H V$ and $R = \RMaps_{\PGLn}(X, \Mn)$.
Note that
\begin{equation} \label{e.invariants} \
k(X)^{\PGLn} = k(\PGLn \times V)^{\PGLn \times H} =
 k(V)^H = k(\alpha_1^{p_1}, \beta_l^{p_1}, 
\dots, \alpha_s^{p_s}, \beta_s^{p_s}) \, . 
\end{equation}
Define elements $\pi_i$ and $\eta_i$ of $R$ by
\begin{equation} \label{e.def-pi-eta}
\begin{aligned}
\pi_i \colon [g,  (\alpha_1, \dots, \alpha_s, \beta_1, \dots, \beta_s)] &
\mapsto   \alpha_i gP_{a_i}g^{-1} \\
\eta_i \colon [g, (\alpha_1, \dots, \alpha_s, \beta_1, \dots, \beta_s)]  
& \mapsto 
\beta_i g D_{\chi_i}g^{-1} \; . 
\end{aligned}
\end{equation}
These elements are well-defined because
$\pi_i(g, v) = \pi_i(gh^{-1}, hv)$ and $\eta_i(g, v) 
= \eta_i(gh^{-1}, hv)$ for every $h \in H$ and $i = 1, \dots, s$;  
see~\eqref{e.pgln}. Note that since $P_{a_i}$ and $D_{\chi_i}$ generate
$\Mn(k)$ as a $k$-algebra, as $i$ ranges from $1$ to $n$  
(cf.~Lemma~\ref{lem.subgr-PGL_n}), 
there exists a dense Zariski dense open subset $X_0 \subset X$ such that
\begin{equation} \label{e.pi-eta}
\text{\rm $\pi_i(x)$ and $\eta_i(x)$ generate $\Mn(k)$ for every $x \in X_0$.}  
\end{equation}
In particular, if $f$ is a central element of $R$ then $f(x)$ is
a scalar matrix for every $x \in X_0$. Consequently, the center $Z(R)$
consists of rational maps $X \brokrarr \Mn$ whose image lies 
in the subspace of scalar matrices.  In other words, 
\begin{equation} \label{e.center}
Z(R) = k(X)^{\PGLn}    
\end{equation}
where, as before, we identify $f \in k(X)^{\PGLn}$
with the $\PGLn$-equivariant rational map $X \brokrarr \Mn(k)$
given by $x \mapsto f(x)I_n$. 

We are now ready to construct an isomorphism between $D_n$ and $R$.
First we identify $D_n$ with the skew-polynomial ring
\[ D_n = Z(R)\{x_1, y_1, \dots, x_s, y_s \} \, ,  \]
where $x_i^{p_i} = \alpha_i^{p_i}$, $y_i^{p_i} = \beta_i^{p_i}$,
$y_ix_i = \zeta_{p_i} x_i y_i$ and all other pair of variables commute.
(Recall that $Z(R)$ is the purely transcendental extension of $k$ 
generated by
$\alpha_1^{p_1}, \beta_1^{p_1}, \dots, \alpha_s^{p_s},\beta_s^{p_s}$;
see~\eqref{e.invariants} and~\eqref{e.center}.) 
Let $\phi \colon D_n \lra R$ be the $Z(R)$-algebra homomorphism 
given by $\phi(x_i) = \pi_i$ and $\phi(y_i) = \eta_i$.
This homomorphism is well-defined because
$\pi_i$ and $\eta_i$ satisfy the same relations as $x_i$ and $y_i$; 
see~\eqref{e.def-pi-eta} and~\eqref{e.pgln}. 

We claim $\phi$ is an isomorphism. Indeed, $\phi$ is injective  
since $D_n$ is a simple algebra. Moreover, since $\dim_k(\Mn) = n^2$,
it is easy to see that $\dim_{Z(R)} \, R \leq n^2$ (see, 
e.g.~, \cite[Lemma 7.4(a)]{tg} for a characteristic-free proof). 
This shows that $\phi$ is an isomorphism and thus completes the proof
of Lemma~\ref{lem.X_D-b}.
\end{pf}

\subsection*{The variety $Y$}

We now define the $\PGL_n$-variety $Y$ by 
 \begin{eqnarray} \label{e.Y}
 Y = \biggl\{ (y_1: \dots: y_m) \in \bbP((\Mn)^m)  \biggm| 
 \begin{array}{l} \sigma^{(i)}(y_1) = \dots = \sigma^{(i)}(y_m) \\  
      \sigma^{(j)}(y_1 \dots y_m) = 0 \end{array}  \biggr\} \, , 
\end{eqnarray}
as in Proposition~\ref{prop:gen.div-alg}. Recall that our
goal is to use Proposition~\ref{prop:gen.div-alg} to show that
the system~\eqref{e.thm3} has no nontrivial solutions.

\begin{lem} \label{lem.thm3}
Under the assumptions of Theorem~\ref{thm3} (i.e., 
$\Char(k)\notdiv n!$, $\sqf(n)\mid m$ and $\sqf(n)\mid i$\/),
$H$ acts on $Y$ without fixed points.
\end{lem}

\begin{pf} The $H$-fixed points in $\bbP((\Mn)^m)$ are of the form 
$y = (y_1: \dots : y_m)$, where each $y_i$ is either 0 
or an element of $\Mn$ which spans a 1-dimensional
character space for $H$. Moreover, the associated characters of all 
non-zero $y_i$ have to be the same. Thus, in view 
of Lemma~\ref{lem.subgr-PGL_n}, there exists an element 
$(a, \chi) \in A \times A^*$ such that $y_i = t_iP_aD_{\chi}$ for some
$t_1, \dots, t_m \in k$.  Note that at least one $t_i$ has 
to be non-zero, since otherwise $y = (0: \dots :0)$ is not 
a well-defined point of $\bbP((\Mn)^m)$.  

Now suppose $y$ is an $H$-fixed point of $Y$. 
Substituting $y_i = t_iP_aD_{\chi}$ into the defining
equations for $Y$, we obtain
\begin{equation} \label{e.z_i}
\left\{ \begin{array}{l}
t_1^i \sigma^{(i)}(P_aD_{\chi}) = \dots = 
t_m^i \sigma^{(i)}(P_aD_{\chi}) \, , \\
t_1\dots t_m \sigma^{(j)}((P_aD_{\chi})^m) = 0 \, . 
\end{array} \right.
\end{equation}
Let $c$ be the order of $(a, \chi)$ in $A \times A^*$. Then
$c\mid\exp(A)$, $\exp(A) = \sqf(n)$, $\sqf(n) \mid m$,
$\sqf(n) \mid i$, and thus, $c\mid m$ and $c\mid i$.
By Lemma~\ref{lem.P_aD_chi}(a), $(P_aD_{\chi})^m=\pm I_n$, and hence,
$\sigma^{(j)}((P_aD_{\chi})^m) \neq 0$.
By Lemma~\ref{lem.P_aD_chi}(c),
$\sigma^{(i)}(P_aD_{\chi}) \neq 0$. Therefore, we can
rewrite~\eqref{e.z_i} as
\[ \left\{ \begin{array}{l}
t_1^i = \dots = t_m^i  \, , \\
t_1\dots t_m = 0 \, . 
\end{array} \right. \]
This system has no solutions other than $t_1 = \dots = t_m = 0$, 
a contradiction.  We conclude that $Y$ has no $H$-fixed points, as claimed. 
\end{pf}

\subsection*{Conclusion of the proof}

In order to complete the proof, it remains to show that 
$X$ has a smooth $H$-fixed point; the desired conclusion will then
follow by applying Proposition~\ref{prop:gen.div-alg} to the 
abelian group $H$ and $\PGLn$-varieties $X$ and $Y$ we introduced above. 

If $D = \UD(n)$ then $X = (\Mn)^2$ (see Lemma~\ref{lem.X_D-a}), and 
the origin is a smooth $H$-fixed point of $X$.

If $D = D_n$ then 
$X = \PGLn \ast_{H} V = (\PGLn \times V)//H$; see Lemma~\ref{lem.X_D-b}.
Since $\PGLn \times V$ is a smooth variely, and $H$ acts freely on it, 
$X$ is also smooth. Moreover, the point of $X$ represented by
$(1, 0) \in \PGLn \times V$, is clearly fixed by $H$. Thus $X$ has a smooth 
$H$-fixed point, as claimed.

This completes the proof of Theorem~\ref{thm3}.
\qed

\subsection*{Refinements}

A slight modification of the above argument proves the following more
general variant of Theorem~\ref{thm3}.

\begin{thm} \label{thm3a}
Let $P(z_1, \dots, z_v) \in k \{ z_1, \dots, z_v \}$ 
be a homogeneous (non-commutative) polynomial of degree $d$
in $v$ variables.  The system of equations
\begin{equation} \label{e.thm3a}
\left\{ \begin{array}{l} \sigma^{(i)}(x_1^u) = \dots = 
\sigma^{(i)}(x_v^u) \\
       \sigma^{(j)}(P(x_1, \dots, x_v)) = 0 \end{array} \right. 
\end{equation}
has no nontrivial solutions in $D_n$ or $\UD(n)$, provided that
 
 \smallskip
(i) $iu$ and $jd$ are divisible by $\sqf(n)$.
 
 \smallskip
(ii) $P(\zeta_1, \dots, \zeta_v) \neq 0$ for any (not necessarily primitive)
$ij$-th roots of unity $\zeta_1, \dots, \zeta_v$.
\end{thm}

Note that if we set $u = 1$, $d= v = m$ and 
$P(z_1, \dots, z_v) = z_1 \dots z_v$, then we recover Theorem~\ref{thm3} 
from Theorem~\ref{thm3a}.

\begin{remark} \label{rem.ref2}
Suppose $K = k(a_1, b_2, \dots, a_l, b_l)$ and 
\[ D = (a_1, b_1)_{r_1} \otimes_K \dots \otimes_K (a_l, b_l)_{r_l} \]
be a tensor product of generic symbol algebras of degree $n = r_1 \dots r_l$. 
Denote the least common multiple of $r_1, \dots, r_l$ by $e$. (Equivalently,
$e$ is the exponent of $D$.) Then the system~\eqref{e.thm3}
has no solutions in $D$ as long as $i$ and $m$ are divisible by $e$.
The proof is the same as above, except that instead of
choosing $H$ and $A$ as in~\eqref{e.H}, we take $H = \phi_A(A \times A^*)$
with $A = (\bbZ/r_1\bbZ) \times \dots \times (\bbZ/r_l \bbZ)$.
Similarly, the system~\eqref{e.thm3a} has no solutions in $D$, provided that
$iu$ and $jd$ is divisible by $e$, and condition (ii) of Theorem~\ref{thm3a}
holds.
\end{remark}

\begin{remark} \label{rem.thm3-fld} Theorem~\ref{thm3} remains true if
$D$ is replaced by the general field extension $L_n/K_n$.
The reason is that there is a natural embedding $\alpha 
\colon L_n \hookrightarrow \UD(n)$ such that 
\[ \alpha \colon \sigma^{(i)}_{L_n/K_n}(y) \mapsto
\sigma^{(i)}_{\UD(n)/Z(n)}(\alpha(y)) \]
for every $y \in L_n$ and every $i = 1, \dots, n$. Indeed, recall
that $\UD(n)$ is generated by two generic $n \times n$-matrices, 
$X = (s_{ij})$ and $Y = (t_{ij})$: we can define
$\alpha(x) = X$ and  $\alpha(a_i)=\sigma^{(i)}(X)$, see, 
e.g.,~\cite[Lemma II.1.4]{procesi1}. If system \eqref{e.thm3}
had a nontrivial solution in $L_n$, it would then 
have a nontrivial solution in $\UD(n)$, contradicting 
Theorem~\ref{thm3}. 
\end{remark}

\begin{remark} \label{rem.prime-to-p} Suppose $\Char(k) = 0$,
$n = p^r$ and $D'$ as
a prime-to-$p$ extension of $D_n$ or $\UD(n)$. Then Theorem~\ref{thm3},
Corollary~\ref{cor4} and Theorem~\ref{thm3a}
remain valid if $D$ is replaced by $D'$. Indeed, let $X$ be as in
Lemma~\ref{lem.X_D-a} (if $D = \UD(n)$) and 
Lemma~\ref{lem.X_D-b} (if $D = D_n$).  Then we can write
$D'$ as $\RMaps_{\PGLn}(X', \Mn)$, where
$X' \brokrarr X$ is a $\PGLn$-invariant rational cover, of degree
prime to $p$. We may assume that $X'$ is smooth and projective. (This 
follows from canonical resolution of singularities; 
see~\cite[Proposition 2.2]{ry2}.)
Since $H$ is a $p$-group, the Going Up Theorem says 
that $X'$ has an $H$-fixed point; see~\cite[Proposition A.4]{ry1}.
The desired conclusion now follows from
Proposition~\ref{prop:gen.div-alg}.
\end{remark}

\section{Proof of Theorem~\ref{thm1}}
\label{sec:pf.thm1}

We may assume without loss of generality that $k$ is an algebraically 
closed field; otherwise we may simply replace $K_n$ and $L_n$ by
$K_n \otimes_k \overline{k}$ and $L_n \otimes_k \overline{k}$ respectively,
where $\overline{k}$ is the algebraic closure of $k$.  

Let $f(x) = x^n + a_1 x^{n-1} + \dots + a_n$ and
$L_n = K_n[x]/(f(x))$, as in~\eqref{general}. 
The normal closure of $L_n$ over
$K_n$ is the field $L'=K_n(x_1,\dots,x_n) = k(x_1,\dots,x_n)$,
where $x_1, \dots, x_n$ are the roots of $f$; they are algebraically
independent over $k$. We will identify $L_n$
with $K_n(x_1)$ by identifying $x \in L_n$ with $x_1 \in k(x_1, \dots, x_n)$
and $a_i$ with $(-1)^i s_i(x_1,\dots,x_n)$, where $s_i$ is the $i$th
elementary symmetric polynomial. 

We shall deduce Theorem~\ref{thm1} as a particular case of
Proposition~\ref{prop:down.Sn}, with $m = 1$, $K = K_n$, $L = L_n$,
$L'$ as above, and $G = \Galois(L'/L_n) = \Sym_n$. 
We will now define the remaining objects that appear in the
statement of Proposition~\ref{prop:down.Sn}, namely the abelian 
subgroup $H$ of $G = \Sym_n$ and the $G$-varieties $X$ and $Y$.

We set $H = H_1 \times H_2$, with
$H_1 = \psi_{A_1}(A_1) \subset \Sym_{n_1}$, 
$H_2 =  \psi_{A_2}(A_2) \subset \Sym_{n_2}$,
as in~\eqref{e.product-group}; here for $i = 1, 2$,
$A_i$ is an abelian subgroup of order $n_i$ and exponent $\sqf(n_i)$.
More precisely, if $n_1 = p_1 \dots p_s$ and $n_2 = q_1 \dots q_t$
are written as products of (not necessarily
distinct) primes then 
\begin{eqnarray} \label{e.primes1-2}
  H_1 \simeq A_1 = (\bbZ/p_1\bbZ) \times \dots \times
(\bbZ/p_s\bbZ)  \nonumber \\  \text{\rm and} \quad \quad \quad \quad \quad  \\
 H_2 \simeq A_2 = (\bbZ/q_1\bbZ) \times \dots \times
(\bbZ/q_t\bbZ)\ .  \nonumber
\end{eqnarray}

We define $X=\bbA^n$, with the natural permutation 
action of $G= \Sym_n$. If we denote the coordinates on $\bbA^n$ by  
$x_1,\dots,x_n$ then $k(X)=k(x_1,\dots,x_n)=L'$ as fields
with $\Sym_n$-action.  The origin is a smooth point of $X$ fixed 
by $\Sym_n$ and, hence, by $H$.

The $S_n$-variety $Y$ is defined as the subvariety 
of $\bbP(\bbA^n)=\bbP^{n-1}$ given by 
\begin{equation} \label{eqn:syst.Y}
\left\{
\begin{aligned}
x_1^{m_1}+\dots+x_n^{m_1} & =0\\
x_1^{m_2}+\dots+x_n^{m_2} & =0\ .
\end{aligned}
\right.
\end{equation}

In order to apply Proposition~\ref{prop:down.Sn}, it is now sufficient
to prove the following:

\begin{lem}
Under the assumptions (i) and (ii) of Theorem~\ref{thm1},
$Y$ has no $H$-fixed points.
\end{lem}

\begin{pf}
By Lemma~\ref{lem.subgr-S_n}, the fixed points $y$ for 
the $H$-action on $\bbP^{n-1} = \bbP^{n_1 + n_2 -1}$ are
of one of the following three types:

\smallskip
Type I:
$y = R_{a, b} =
(\underbrace{a: \dots: a}_{\text{$n_1$ times}}:
\underbrace{b: \dots : b}_{\text{$n_2$ times}})$,
for some $a, b \in k$, not both 0.

\smallskip
Type II: $y = (R_{\chi}, 0)  = 
(\chi(\alpha_1): \dots: \chi(\alpha_{n_1}): 0: \dots: 0)$, 
where $H_1 =
\{ \alpha_1, \dots, \alpha_{n_1} \}$ 
and  $\chi$ is a character of $H_1$.

\smallskip
Type III: $y = (0, R_{\eta}) = 
(0: \dots: 0: \eta(\beta_1): \dots: \eta(\beta_{n_2}))$, 
where $H_2 =
\{ \beta_1, \dots, \beta_{n_2} \}$ and  $\eta$ is a character of $H_2$.

\smallskip
\noindent

Consider a point of type I. 
Substituting the coordinates of $R_{a, b}$ into \eqref{eqn:syst.Y},
we see that $R_{a, b}$ lies in $Y$
if and only if $(a, b)$ is a nontrivial solution of
the homogeneous system
\begin{equation} \label{sys.tr}
\left\{ \begin{array}{l} n_1 a^{m_1} + n_2 b^{m_1} = 0 \\
n_1 a^{m_2} + n_2 b^{m_2} = 0 \, . \end{array} \right. 
\end{equation}
An elementary computation shows that under assumption (i) of 
Theorem~\ref{thm1} this system has no nontrivial solutions. Hence
we conclude that no point of type I can lie on $Y$.

We now turn to points of types II and III. Since
$H_1$ has exponent $\sqf(n_1)$, we see that
$\chi(\alpha_i)^{\sqf(n_1)} = 1$ for every $\alpha_i \in H_1$.
It follows from the assumptions of Theorem~\ref{thm1} that $n_1\ne0$
in $k$ and either $m_1$ or $m_2$ is divisible by $\sqf(n_1)$;
consequently, $(R_{\chi}, 0)$ does not lie on $Y$.
Similarly, $(0, R_{\eta})$ does not lie on $Y$.
Hence, no point of type II or III lies on $Y$.
This completes the proof of the lemma and thus of Theorem~\ref{thm1}.
\end{pf}

\begin{remark} \label{rem.thm1-div} Theorem~\ref{thm1} fails if the 
field extension $L_n/K_n$ is replaced by the generic division algebra
$\UD(n)$. Suppose, for simplicity, that $k$ is an algebraically closed 
field of characteristic zero. Then, by a theorem of Wedderburn,
$\UD(3)$ is cyclic; thus it has an
elements $x$ and $y$ such that $x = \zeta_3 yxy^{-1}$, where
$\zeta_3$ is a primitive cube root of 1. It is now easy to see that
$\tr(x) = \tr(x^2) = 0$. On the other hand, Theorem~\ref{thm1} with
$n_1 = m_1 = 1$ and $n_2 = m_2 = 2$, says that no such element can 
exist in $L_3$.  

Another example of this kind can be constructed for $n = 6$.
The algebra $D = \UD(6)$ is known to be cyclic; 
hence, it has a non-zero element $z$ such that $\tr(z^i) = 0$
for $i = 1, \dots, 5$. On the other hand, Theorem~\ref{thm1} says 
that the systems $\tr(x) = \tr(x^5) = 0$ or $\tr(x^2) = \tr(x^4)=0$
have no solutions in $L_6^*$.
\end{remark}

\begin{remark} \label{rem.etale-alg} Let $n_1 = p_1 \dots p_s$ and
$n_2 = q_1 \dots q_t$, where $p_1, \dots, p_s, q_1, \dots, q_t$ are
(not necessarily distinct) primes. Suppose $z_1, \dots, z_s$ and
$w_1, \dots, w_t$ are independent variables over $k$. Set
$E_1 = k(z_1, \dots, z_s, w_1^{q_1}, \dots, w_t^{q_t})$,
$E_2 = k(z_1^{p_1}, \dots, z_s^{p_s}, w_1, \dots, w_t)$, and
$F = k(z_1^{p_1}, \dots, z_s^{p_s}, w_1^{q_1}, \dots, w_t^{q_t})$. 
Then we can replace $L_n/K_n$ by the $n$-dimensional
etale $F$-algebra $E = E_1 \oplus E_2$ (cf.~\cite[Section~4]{reichstein})
in the statement of~Theorem~\ref{thm1}. In other words,

\smallskip{\narrower\noindent
{\em under assumptions (i) and (ii) of Theorem~\ref{thm1} the system
of equations $\tr(x^{m_1}) = \tr(x^{m_2}) = 0$ has no nontrivial
solutions in $E$.}

\smallskip}

The role played by $E$ in this setting
is analogous to the role played 
by $D_n$ in the setting of Theorem~\ref{thm3}. 
In particular, one can show that $E = \RMaps_{\Sym_n}(X, \bbA^n)$, where
$X = \Sym_n \ast_H V$, $V$ is a faithful
$(s + t)$-dimensional linear representation of $H = H_1 \times H_2$, 
and the algebra structure on $\RMaps_{\Sym_n}(X, \bbA^n)$ is 
induced from the algebra structure on
$\bbA^n = \underbrace{k \oplus \dots \oplus k}_{\text{$n$ times}}$
(compare with Lemma~\ref{lem.X_D-b}). 
Since $X$ has a smooth $H$-fixed point (namely, the point
represented by $(id, 0) \in \Sym_n \times V$), the rest of 
our argument goes through unchanged. 
\end{remark}
  
\section{Systems of the form $\sigma^{(m_1)}(x) = 
\sigma^{(m_2)}(x) = 0$}
\label{sect5}

We do not know whether or not the system $\tr(x^{m_1}) = \tr(x^{m_2}) = 0$
may be replaced by the system 
\begin{equation} \label{e.prop1a}
\sigma^{(m_1)}(x) = \sigma^{(m_2)}(x) = 0 \, . 
\end{equation}
in the statement of Theorem~\ref{thm1}.  (Such a result would be of 
interest, since it would mean that the general polynomial of degree $n$ 
cannot be transformed, by a Tschirnhaus substitution, into a polynomial 
$t^n + b_1t^{n-1} + \dots + b_n$, with $b_{m_1} = b_{m_2} = 0$.) Every 
step of our proof of Theorem~\ref{thm1} goes through in this case, except 
that the system~\eqref{sys.tr} is replaced by the system
\begin{equation} \label{sys.sigma}
\left\{ \begin{array}{l} 
  s_{m_1}(a, \dots, a, b, \dots, b) = 0 \\
  s_{m_2}(a, \dots, a, b, \dots, b) = 0 \, , \end{array} \right. 
\end{equation}
where $(a, \dots, a, b, \dots, b)$ stands for 
$(\underbrace{a, \dots, a}_{\text{$n_1$ times}},
\underbrace{b, \dots, b}_{\text{$n_2$ times}})$  and
$s_i$ denotes the $i$th elementary symmetric polynomial. Thus: 

\begin{prop} \label{prop1a} Let $n_1$ and $n_2$ be positive integers
prime to $\Char(k)$,
and $L_n/K_n$ be the general field extension of degree $n= n_1 + n_2$.
Then the system~\eqref{prop1a} 
has no nontrivial solutions $x \in L_n^*$, provided that 
each $\sqf(n_i)$ ($i = 1, 2$) divides $m_1$ or $m_2$ and 
the system~\eqref{sys.sigma} has no nontrivial solutions 
$(a, b) \in k^2$.
\end{prop}

Of course, this result is less satisfying than Theorem~\ref{thm1} because
we do not know for what values of $n_1$, $m_1$, 
$n_2$ and $m_2$ the system~\eqref{sys.sigma} has no nontrivial 
solutions. (The analogous question for the 
system~\eqref{sys.tr} is quite easy: the answer is given by
condition (i) of Theorem~\ref{thm1}.) Nevertheless, 
for low values of $n$, Proposition~\ref{prop1a} gives us a rather 
complete picture. We shall give two such examples below.

Before preceeding with the examples, we record a simple observation.

\begin{remark} \label{rem.inverse}
Let $E/F$ be a field extension of degree $n$. Multiplying
\eqref{ci} by $\det((\lambda x)^{-1})$, we easily obtain the identity
$\sigma^{(n-i)}(x^{-1}) = \sigma^{(i)}(x)/\sigma^{(n)}(x)$.
In particular, if $x \in E$ satisfies~\eqref{e.prop1a} 
then $\sigma^{(n-m_1)}(x^{-1}) = \sigma^{(n-m_2)}(x^{-1}) = 0$. 
\qed
\end{remark}

\begin{example} \label{ex.deg5}
Let $L_5/K_5$ be the general field extension of degree $5$ and let
$1\leq m_1 < m_2 \leq 5$. Then the system~\eqref{e.prop1a}
has a nontrivial solution $x \in L_5^*$ if and only if $(m_1, m_2) = (1, 3)$
or $(2, 4)$.
\end{example}

\begin{pf} By the theorem of Hermite cited in Example~\ref{ex1.1},
the system~\eqref{e.prop1a} has a solution $0 \neq x \in L_5$ for 
$(m_1, m_2) = (1, 3)$. Then $x^{-1}$ is a solution to~\eqref{e.prop1a}
with $(m_1, m_2) = (2, 4)$; see Remark~\ref{rem.inverse}.

It remains to show that there are no solutions for any other values of 
$m_1$ and $m_2$. Indeed, we may assume without loss of generality that
$m_2 \neq 5$, since $\sigma^{(5)}(x) = - \det(x) \neq 0$ for any $x 
\in L_5^*$.  The remaining 
possibilities for $(m_1, m_2)$ are: $(1, 2)$, $(1, 4)$, $(2, 3)$, 
and $(3, 4)$. In view of Remark~\ref{rem.inverse}, we only need to
consider $(1, 2)$, $(1, 4)$ and $(2, 3)$.

\smallskip
$(m_1, m_2) = (1,2)$. By Newton's formulas the system
$\sigma^{(1)}(x) = \sigma^{(2)}(x) = 0$
is equivalent to $\tr(x) = \tr(x^2) = 0$. The latter system has no solutions
by Theorem~\ref{thm1} with $n_1 = 1$ and $n_2 = 4$. (Alternatively, use
Proposition~\ref{prop1a} with $n_1 = 1$, $n_2 = 4$ or
appeal to~\cite[Theorem 1.3(b)]{reichstein}, with $p = 2$ and $m = 2$.)

\smallskip
$(m_1, m_2) = (1, 4)$. 
Apply Proposition~\ref{prop1a} with $n_1 = 1$ 
and $n_2 = 4$. In this case~\eqref{sys.sigma} reduces to 
\[  \left\{ \begin{array}{l} 
 s_1(a, b, b, b, b) = a + 4b = 0 \\
 s_4(a, b, b, b, b) = b^4 + 4ab^3 = 0 \, . \end{array} \right. \] 
It is easy to see that this system has no nontrivial solutions. 
(Alternatively, use~\cite[Theorem 6.1b]{reichstein}.) 

\smallskip
$(m_1, m_2) = (2, 3)$. Apply Proposition~\ref{prop1a} with $n_1 = 2$ 
and $n_2 = 3$. In this case~\eqref{sys.sigma} becomes
\[  \left\{ \begin{array}{l} 
 s_2(a, a, b, b, b) = a^2 + 6ab +3b^2 = 0 \\
 s_3(a, a, b, b, b) = 3a^2b + 6ab^2 + b^3= 0 \, . \end{array} \right. 
\]
This system has no nontrivial solutions.
\end{pf}
 
\begin{example} \label{ex.deg6}
Let $L_6/K_6$ be the general field extension of degree $6$ and let
$1\leq m_1 < m_2 \leq 6$. Then the system~\eqref{e.prop1a}
has a nontrivial solution $x \in L_5^*$ if and only if $(m_1, m_2) = (1, 3)$
or $(3, 5)$.
\end{example}

\begin{pf} The existence of solutions for $(m_1, m_2) = (1, 3)$ and
$(3, 5)$ follows from Example~\ref{ex1.1} and Remark~\ref{rem.inverse}.

We may assume $m_2 \leq 5$ because $\sigma^{(6)}(x) = \det(x) \neq 0$ for
any $x \in L_6^*$.
It is now enough to show that there are no solutions for $(m_1, m_2) =
(1, 2)$, $(1, 4)$, $(1, 5)$, $(2, 3)$, and $(2, 4)$; the remaining cases
follow from these by Remark~\ref{rem.inverse}. 

\smallskip
$(m_1, m_2) = (1, 2)$. In this case~\eqref{sys.sigma} is
equivalent to $\tr(x) = \tr(x^2) = 0$. The latter system has 
no solutions by Theorem~\ref{thm1} with $n_1 = 2$ and $n_2 = 4$.  
(Alternatively, use Proposition~\ref{prop1a} with $n_1 = 1$, $n_2 = 4$ or
appeal to~\cite[Theorem 1.3(c)]{reichstein}, with $p =2$, $m = 2$ 
and $l = 1$.)

\smallskip
$(m_1, m_2) = (1, 4)$. Apply Proposition~\ref{prop1a} with
$n_1 = 2$, $n_2 = 4$. In this case~\eqref{sys.sigma} 
reduces to $2a+4b = 6a^2 b^2 + 8ab^3 + b^4 = 0$. This system has 
no nontrivial solutions. 

\smallskip
$(m_1, m_2) = (1, 5)$.  Apply Proposition~\ref{prop1a} with
$n_1 = 1$, $n_2 = 5$. In this case~\eqref{sys.sigma} 
reduces to
$a+5b = 5ab^4 + b^5 = 0$.
There are no nontrivial 
solutions.  (Alternatively, use~\cite[Theorem 1.3(b)]{reichstein} with
$p = 5$.)

\smallskip
$(m_1, m_2) = (2, 3)$.  Apply Proposition~\ref{prop1a} with
$n_1 = 2$, $n_2 = 4$. In this case~\eqref{sys.sigma} becomes
$a^2 + 8ab + 6b^2 = 4a^2b + 12ab^2 + 4b^3 = 0$. There
are no nontrivial solutions.

\smallskip
$(m_1, m_2) = (2, 4)$.  Use Proposition~\ref{prop1a} with
$n_1 = 2$, $n_2 = 4$. In this case~\eqref{sys.sigma} becomes
$a^2 + 8ab + 6b^2 = 6a^2b^2 + 8ab^3 + b^4 = 0$. Once again, there
are no nontrivial solutions.
\end{pf}

\section{A further generalization}
\label{sect6}

In this section we will show that the assumption that the $G$-variety $Y$
in Proposition~\ref{prop:down.Sn} has no fixed points can sometimes
be weakened. We will present a general result extending
Proposition~\ref{prop:down.Sn} and illustrate it with an example.
One can generalize Proposition~\ref{prop:gen.div-alg} in
a similar manner; we leave the details to an interested reader.

In this section we assume that $k$ is algebraically closed.

\begin{prop} \label{prop:down.Sn.gen'lized}
Assume

\smallskip
(i) $L/K$ is a separable field extension of degree $n$,
$L'$ is the normal closure of $L$ over $K$, $G = \Gal(L', K)$,
and $H$ is an abelian subgroup of $G$,

\smallskip
(ii) $Y\supset Z$ are subvarieties of $(\bbA^n)^m$
given, respectively, by systems of $G$-invariant polynomial equations
$P_1=\dots=P_s=0$ and $Q_1=\dots=Q_r=0$,  

\smallskip
(iii) there exists a complete $H$-variety $W$ without $H$-fixed 
points and a regular $H$-equivariant map $h \colon Y-Z \to W$, 
and

\smallskip
(iv) there exists a $G$-variety $X$ such that such that $k(X)=L'$ as 
fields with $G$-action, and $X$ has  a smooth $H$-fixed point.

\smallskip
Then any solution $(a_1,\dots,a_m)\in L^m$ of the system
\begin{equation} \label{eqn:gen.in.field.P.gen'lized}
P_1(x_1,\dots, x_m)=\dots=P_s(x_1,\dots,x_m)=0
\end{equation}
also satisfies  the system
\begin{equation} \label{eqn:gen.in.field.Q.gen'lized}
Q_1(x_1,\dots,x_m)=\dots=Q_r(x_1,\dots,x_m)=0\ .
\end{equation}
\end{prop}

Note that since $Z \subset Y$, the ideal
$(Q_1,\dots,Q_r)\subset k[(\bbA^n)^m]$ contains some power of the
ideal $(P_1,\dots,P_s)$. Hence, any solution of
\eqref{eqn:gen.in.field.Q.gen'lized} in $L^n$ is a solution of
\eqref{eqn:gen.in.field.P.gen'lized}.
Proposition~\ref{prop:down.Sn.gen'lized} asserts that
under assumptions (i)--(iv), the opposite is also true. 

\begin{pf}
Given a solution $(a_1,\dots,a_m)$ of
\eqref{eqn:gen.in.field.P.gen'lized}, we construct
a rational map $f\colon X\brokrarr Y \subset (\bbA^n)^m$, as in 
the proof of Proposition~\ref{prop:down.Sn}.
If $(a_1,\dots,a_m)$ does not satisfy
\eqref{eqn:gen.in.field.Q.gen'lized}, then $f(X)\not\subset Z$ and
hence, the composition
$X\overset{f}{\brokrarr}Y\overset{h}{\brokrarr}W$ is a well-defined
$H$-equivariant rational map.
As $X$ has a smooth $H$-fixed point, Theorem~\ref{going-down} 
says that $W$ also has one, a contradiction.
\end{pf}

\begin{remark} \label{ex:back.to.hom}
To see that Proposition~\ref{prop:down.Sn} is a special 
case of Proposition~\ref{prop:down.Sn.gen'lized}, assume that the 
polynomials $P_1,\dots,P_s$ are homogeneous, so that $Y$ is a cone in
$(\bbA^n)^m$, and $Z$ is the origin in $(\bbA^n)^m$. Note that
the origin of $(\bbA^n)^m$ can be cut out by $G$-invariant homogeneous 
polynomials (this is true for any finite group representation), thus
we can choose $Q_1, \dots, Q_r \in k[(\bbA^n)^m]^G$ to be generators of
the ideal of the origin in $k[(\bbA^n)^m]$.

Let $W \subset\bbP((\bbA^n)^m)$ be the projectivisation of the cone
$Y$, and $h\colon Y - Z \lra W$ the natural projection.
If $W$ has no $H$-fixed points, and $X$ has a smooth $H$-fixed point
then Proposition~\ref{prop:down.Sn.gen'lized} implies that the system
\eqref{eqn:gen.in.field.P.gen'lized} has no solutions, except for
$x_1=\dots=x_m=0$.  This is precisely the statement of
Proposition~\ref{prop:down.Sn}.
\end{remark}

\begin{remark} \label{ex:gen.resol}
Proposition~\ref{prop:down.Sn.gen'lized} can be
applied in the following situation.
Suppose that $Z$ is the singular set of $Y$.
Let $\tilde Y$ be the closure of $Y\subset(\bbA^n)^m=\bbA^{nm}$ in
$\bbP^{nm}\supset\bbA^{nm}$; note that the $G$-action on $(\bbA^n)^m$
extends to a regular $G$-action on $\bbP^{nm}$, and $\tilde Y$ is
$G$-invariant.
Let $\pi \colon W \lra\tilde Y$ be the canonical resolution of 
singularities. Such a resolution is known to exist if $\Char(k)=0$; 
see the discussion and the references in~\cite[Section 3]{ry1}.
Note that $\pi$ is an isomorphism over $Y - Z$ and thus
we can take $h=\pi^{-1} \colon Y - Z \lra W$.
If $W$ has no 
$H$-fixed points then Proposition~\ref{prop:down.Sn.gen'lized} applies.
\end{remark}

\begin{example} \label{example:high.powers}
{\em Suppose $n$ is prime and $n \ne \Char(k)$.
Then for any $c\in k$ the equation 
\begin{equation} \label{eqn:high.powers}
\sum_{i = 1}^{n-1} \sigma^{(i)}(x)^n \sigma^{(n)}(x)^{n-1-i} + 
c\sigma^{(n)}(x)^{2n-2}  = 0 
\, , \end{equation}
has no nontrivial solutions in the general field extension $L_n/K_n$;
see~\eqref{general}. Here $\sigma^{(i)}$ stands 
for $\sigma_{L_n/K_n}^{(i)}$.}
\end{example}

\begin{pf} We may assume without loss of generality that $k$ is
algebraically closed, and thus, contains the roots of unity.

First consider the case $c\ne0$.
We apply Proposition~\ref{prop:down.Sn.gen'lized} in the following
setting: $K = K_n$, $L = L_n$,  $G = \Sym_n$, $X = \bbA^n$ with 
the natural $\Sym_n$-action, $H$ = the cyclic subgroup of $\Sym_n$
generated by the $n$-cycle $h=(1\,2\,\dots\,n)$, $s=m=1$, and
$P_1=
s_1^n s_n^{n-2}+s_2^n s_n^{n-3}+\dots+s_{n-1}^n+cs_n^{2n-2}$,
where $s_i$ denotes the $i$th elementary symmetric polynomial in the
coordinates $x_1,\dots,x_n$ in $\bbA^n$. (To construct $P_1$, we 
replaced $\sigma^{(i)}(x)$ by $(-1)^i s_i(x_1, \dots, x_n)$ in the
left hand side of~\ref{eqn:high.powers}.)
Note that $P_1$ is not homogeneous in $x_1, \dots, x_n$ as $c\ne0$.

We take $Z$ to be the origin in $\bbA^n$.
Similarly to Remark~\ref{ex:gen.resol}, let $\tilde Y$ the closure of
$Y\subset A^n$ in $\bbP^n$; then the $H$-action on $\tilde Y-Z$ is
free.
Let $\widetilde{\bbP^n}\to\bbP^n$ be the blowup of $Z$; we identify
its exceptional divisor $S$ with $\bbP^{n-1}$.
Let $Y'$ be the strict transform of $\tilde Y$; then
$Y'\to\tilde Y$ is a blowup centered at $Z$, and $S\cap Y'$
is the hypersurface in $\bbP^{n-1}$ given
by the homogeneous equation $\overline P_1=0$ where
$\overline P_1=s_1^n s_n^{n-2}+s_2^n s_n^{n-3}+\dots+s_{n-1}^n$
is the initial form of $P_1$.

The intersection $S\cap Y'$ contains $H$-fixed points
$q_\zeta=(1:\zeta:\zeta^2:\dots:\zeta^{n-1})$ for each $n$-th root of
unity $\zeta\ne1$.  Let $W\to Y'$ be the blowup 
of these $n-1$ points. We claim that $W$ has no $H$-fixed points.

To see this, consider the hypersurfaces $S_i\subset\widetilde{\bbP^n}$
for $i=1,\dots,n-1$ which are the closures in $\widetilde{\bbP^n}$ of
the hypersurfaces in $\bbA^n-Z$ given by the equations $s_i=0$. 
For each $i$, the intersection $S_i\cap S$ is the hypersurface in
$S=\bbP^{n-1}$ given by the homogeneous equation $s_i=0$; in particular,
each $S_i$ passes through $q_\zeta$. 
Consider the $(n-1)\times(n-1)$ Jacobian determinants
$D_l(q_\zeta)=\det(\partial s_i/\partial x_j)(q_\zeta)$,
where $i=1,\dots,n-1$ and $j=1,\dots,\widehat l,\dots,n$.
By Newton's formulas
$D_l(q_\zeta)=\det(\partial p_i/\partial x_j)(q_\zeta)$, where
$p_i=x_1^i+\dots+x_n^i$. The latter determinant is a Vandermonde
determinant, which does not vanish at $q_\zeta$. This shows that
the hypersurfaces $S_i\cap S$ are smooth and
intersect transversely (in $S = \bbP^{n-1}$) at each $q_\zeta$;
hence $S_1, \dots, S_{n-1}$ and $S$ are smooth and intersect
transversely (in $\widetilde{\bbP^n}$) at each $q_{\zeta}$.

Thus the tangent spaces
$T_{q_\zeta}(S_1), \dots, T_{q_\zeta}(S_{n-1})$, 
together with $T_{q_\zeta}(S)$, form a system of coordinate
hyperplanes in $T_{q_\zeta}(\widetilde{\bbP^n})$. Since each $S_i$
is $H$-invariant, the linear $H$-action on
$T_{q_\zeta}(\widetilde{\bbP^n})$ is diagonalized in this coordinate system.
The group $H$ acts by different characters on each of the coordinate
directions; in fact, $h$ acts by multiplication by $\zeta^i$ on
$T_{q_\zeta}(\widetilde{\bbP^n})/T_{q_\zeta}(S_i)$, and trivially on
$T_{q_\zeta}(\widetilde{\bbP^n})/T_{q_\zeta}(S)$.
Identifying the exceptional divisor $E_{q_{\zeta}}$
of the blowup of $\widetilde{\bbP^n}$ centered at $q_\zeta$, with
$\bbP(T_{q_\zeta}(\widetilde{\bbP^n}))$, we see
that the $H$-fixed points on $E_{q_{\zeta}}$
are the points of $\bbP(T_{q_\zeta}(\widetilde{\bbP^n}))$ that
correspond to the directions of the coordinate axes in
$T_{q_\zeta}(\widetilde{\bbP^n})$.
The exceptional divisor of $W$ over $q_\zeta$ is the
projectivisation of the tangent cone to $Y'$ at $q_\zeta$, and
the latter does not contain the coordinate axes. We conclude
that $W$ does not have $H$-fixed points, as claimed.

Thus, we may apply Proposition~\ref{prop:down.Sn.gen'lized};
it shows that the equation \eqref{eqn:high.powers} has no nontrivial
solutions, similarly to Remark~\ref{ex:back.to.hom}.

In case $c=0$, we need to make the following changes.
Now $Y$ is an affine cone; we take $Z$ to be
the union of $(n-1)!$ lines that correspond to the points
$(\zeta_1:\dots:\zeta_n)\in\bbP^{n-1}$ where $\zeta_1,\dots,\zeta_n$
are different $n$th roots of unity; this includes the
lines that correspond to the points $q_\zeta$.
Now let $Y'$ be the blowup of $\tilde Y$ at the origin as
before, and $W$ be the blowup of $Y'$ at the lines
that make up the strict transform of $Z$ in $Y'$.
(Alternatively, we may take the route similar to
Remark~\ref{ex:back.to.hom} and set $W$ to be the blowup of $\bbP(Y)$
at the points $q_\zeta$.)
Then $W$ does not have $H$-fixed points, and
Proposition~\ref{prop:down.Sn.gen'lized}
shows that any $x \in L_n$ satisfying~\eqref{eqn:high.powers} also
satisfies the system \eqref{eqn:gen.in.field.Q.gen'lized}, which in our
case is 
\begin{equation} \label{e.system2}
\sigma^{(1)}(x)=\dots=\sigma^{(n-1)}(x)=0 \; .
\end{equation}
One can now show directly that $L_n$ does not have a non-zero element $x$
satisfying~\eqref{e.system2}; otherwise $L_n/K_n$ would have to
be a cyclic extension, a contradiction. 
Alternatively, one can show that the system~\eqref{e.system2}
has no nontrivial solutions by applying
Proposition~\ref{prop:down.Sn.gen'lized} one more time, as follows:
\begin{itemize}
\item[---] take the new $H$ to be any cyclic subgroup of $G=\Sym_n$ of
order different from $n$ and $1$;
\item[---] the new $Y$ to be the old $Z$, i.e., $P_i 
= s_i(x_1, \dots, x_n)$
for $i = 1, \dots, n-1$.

\item[---] the new $Z$ to be the origin in $\bbA^n$, i.e., 
$Q_j = s_j(x_1, \dots, x_n)$ for $j = 1, \dots, n$.

\item[---] the new $W$ to be the normalization of $Z$, i.e., the disjoint
union of $(n-1)!$ lines.
\end{itemize}
Applying Proposition~\ref{prop:down.Sn.gen'lized} we see that
the system~\eqref{e.system2} has no nontrivial solutions and, hence,
neither does equation~\eqref{eqn:high.powers}.
\end{pf}


\section{Equations in octonion algebras}
\label{sect7}

\subsection*{Preliminaries}
Let $F$ be a field of characteristic $\neq 2$.
Recall that for any $0 \neq a, b, c \in F$, 
the octonion (or Cayley---Dickson) algebra $\Or_F(a, b, c)$ is defined 
as follows.  The quaternion algebra
\[  (a, b)_2 = F\{i, j \}/(i^2 = a, j^2 = b, ji = -ij) \]
is equipped with an involution $x \rightarrow \overline{x}$ 
given by
\begin{equation} \label{eqinv}
\overline{x_0 + x_1i + x_2 j + x_3ij} = x_0 -x_1i - x_2 j -x_3 ij 
\end{equation}
for any $x_0, \ldots, x_3 \in F$.
Now $\Or_F(a, b, c) \stackrel{{\rm def}}{=} (a, b)_2 \oplus (a, b)_2 l$ is an
8-dimensional $F$-algebra with (non-associative) multiplication  
given by 
$(x+yl)(z+wl) = (xz + c\overline{w}y) + (wx + y\overline{z})l$.
The involution (\ref{eqinv}) extends from $(a, b)_2$ to $\Or_F(a, b, c)$
via $\overline{x + yl} = \overline{x} - yl$. The algebra $\Or_F(a, b, c)$
is also equipped with $F$-valued trace and norm functions given by
$\tr(x) = x + \overline{x}$ and $n(x) = x \overline{x} = 
\overline{x} x$ such that $ x^2 - \tr(x)x + n(x) = 0 $ 
for any $x \in \Or_F(a, b, c)$; we can think of $\tr(x)$ as 
$\sigma^{(1)}(x)$ and $n(x)$ as $\sigma^{(2)}(x)$. 
Note that $\tr(x)$ is intrinsically defined in $\Or_K(a, b, c)$, 
i.e., $\tr(x) = \tr(\sigma(x))$, where $\sigma$ is a $K$-algebra 
automorphism in $\Or_K(a, b, c)$; the same is true of $n(x)$. 
For a more detailed description of octonion algebras
we refer the reader to~\cite{schaefer}.

Two octonion algebras will be of particular interest to us:
the {\em split} algebra $\Or_F(1, 1, 1)$ over $F$ and
the {\em generic} algebra $\Or_{gen} = \Or_K(a, b, c)$, 
where $K = k(a, b, c)$ 
and $a, b, c$ are algebraically independent over $k$.

By a theorem of Zorn~\cite[III.3.17]{schaefer}, 
any $8$-dimensional $F$-algebra $A$ such that
$A \otimes_F F' \simeq \Or_F(1, 1, 1)$ for some field extension $F'/F$,
is necessarily isomorphic to $\Or_F(a, b, c)$ for some $a, b, c \in F^*$.
This means that octonion algebras are ``forms" of the split octonion algebra
$\Or_k(1, 1, 1)$ in the same way as central simple algebras are ``forms"
of the matrix algebra $\Mn(k)$.

\subsection*{$G_2$-equivariant maps}
 From now on we shall assume the base field $k$ to be algebraically 
closed and of characteristic $\neq 2$.

Recall that the automorphism group of the split octonion algebra
$\Or =\Or_k(1, 1, 1)$ is the exceptional group $G_2$. Octonion 
algebras are related to $G_2$-varieties in the same way as 
central simple algebras are related to $\PGLn$-varieties.
In particular, if $k$ is
of characteristic 0 then any octonion
algebra whose center is a finitely generated field extension of $k$
can be written in the form $\RMaps_{G_2}(X, \Or)$, where
$\Or$ is viewed as an 8-dimensional vector space with the natural
$G_2$-action and $X$ is a generically free $G_2$-variety, uniquely 
determined up to birational isomorphism.

 From now on, let $H \simeq (\bbZ/2)^3$ be the subgroup of $G_2$
generated by $\tau_1$, $\tau_2$ and $\tau_3$, where
\begin{equation} \label{e.oct1}
\begin{array}{lll}
\tau_1(i) = -i \, , & \tau_1(j) = j \, , & \tau_1(l) = l \, ; \\
\tau_2(i) = i \, , & \tau_2(j) = -j \, , & \tau_2(l) = l \, ; \\
\tau_3(i) = i \, , & \tau_3(j) = j \, , & \tau_3(l) = -l \, . 
\end{array} 
\end{equation}

\begin{lem} \label{lem.oct2}  The generic octonion algebra
$\Or_{gen}$ is isomorphic to $\RMaps_{G_2}(X, V)$, where
$X = G_2 \ast_H V$ and $V = \Span \{ i, j, k \}$ is the 3-dimensional 
faithful representation of $H$ given by~\eqref{e.oct1}. 
\end{lem}

\begin{pf} The proof is similar to the proof of Lemma~\ref{lem.X_D-b},
so we will only outline it below.

Let $\alpha, \beta, \gamma$ be the coordinates
of $V$ relative to the basis $\{i, j, l \}$,
let $R = \RMaps_{G_2}(X, \Or)$ and let
$\pi_1,\pi_2,\pi_3 \colon X \brokrarr \Or$ be
the elements of $R$ given by 
\begin{equation} \label{e.oct2} \begin{array}{l}
\pi_1 \colon [g, (\alpha, \beta, \gamma)] \mapsto \alpha g(i) \\
\pi_2 \colon [g, (\alpha, \beta, \gamma)] \mapsto \beta g(j) \\
\pi_3 \colon [g, (\alpha, \beta, \gamma)] \mapsto \gamma g(l) \, . 
\end{array} \end{equation}
It is easy to see that these maps are well-defined, i.e.
$\pi_a(g, v) = \pi_a(gh^{-1}, hv)$. 
Let
\[  K = k(X)^{G_2} = k(V)^H = k(\alpha^2, \beta^2, \gamma^2) \; . \]
We now identify $\Or_{gen}$ with 
$\Or_K(\alpha^2, \beta^2, \gamma^2)$,
and define $\phi \colon \Or_{gen} \lra R$ by $\phi(i) = \pi_1$,
$\phi(j) = \pi_2$ and $\phi(l) = \pi_3$. Then $\phi$ is well-defined;
see~\eqref{e.oct2}. Since $\Or$ is a (non-associative) division
algebra, $\phi$ is injective. To see that $\phi$ is an isomorphism, 
we only need to show that $\dim_K(R) \leq 8$; this follows
from~\cite[Lemma 7.4(a)]{tg}.
\end{pf}

\subsection*{$G_2$-invariant polynomials}
Consider the diagonal $G_2$-action on the $8m$-dimensional $k$-vector 
space $W = \Or^m$. Let $P \in k[W]^{G_2}$ be a $G$-invariant
polynomial and let $A = \Or_F(a, b, c)$ be an octonion algebra. 
Identifying $A$ with an $F$-subalgebra of 
$A \otimes_F F' \simeq \Or_{F'}(1, 1, 1)$, where
$F' = F (\sqrt{a}, \sqrt{b}, \sqrt{c})$, we can define
$P(a_1, \dots, a_m)$ for any $a_1, \dots, a_m \in A$. 
Arguing as in Lemma~\ref{lem:inv.polyn}, we see that
$P(a_1, \dots, a_m)$ is well-defined and lies in $F$ for any  
$a_1, \dots, a_m \in A$. (This also follows from a theorem of
Schwarz~\cite[(3.23)]{schwarz},
which asserts that $k[W]^{G_2}$ is generated by elements of 
the form $\tr(M)$, where $M$ is a monomial in $u_1, \dots, u_m 
\in \Or$.)

\begin{prop} \label{prop.oct3} 
Let $H \simeq (\bbZ/2)^3$ be the subgroup
of $G_2$ defined in~\eqref{e.oct1}.  Suppose the subvariety $Y$ 
of $\bbP(\Or^m)$, cut out by homogeneous $G_2$-invariant 
polynomials $P_1 = \dots =  P_r = 0$, does not have an $H$-fixed point.
Then the system 
\begin{equation} \label{e.oct-system}
 P_1(x_1, \dots, x_m) = \dots = P_r(x_1, \dots, x_m) = 0
\end{equation}
has no non-trivial solutions in any octonion algebra of the form
$\RMaps_{G_2}(X, \Or)$, where $X$ is a $G_2$-variety with a 
smooth $H$-fixed point. In particular, the system~\eqref{e.oct-system}
has no nontrivial solutions in the generic octonion algebra $O_{gen}$. 
\end{prop}

\begin{pf} We argue as in the proof of Proposition~\ref{prop:gen.div-alg}.
Assume, to the contrary, that $(a_1, \dots, a_m)$ 
is a nontrivial solution of~\eqref{e.oct-system}.  Each $a_i$ 
is a $G_2$-equivariant rational map $X \brokrarr \Or^m$; 
together they define a $G_2$-equivariant
rational map $f \colon X \brokrarr Y 
\subset \bbP(\Or^m)$. Applying the Going Down Theorem~\ref{going-down},
we obtain a contradiction. 

This proves the first assertion of the proposition. The second assertion  
follows from Lemma~\ref{lem.oct2}. Indeed, the variety $X = G_2 \ast_H V$ 
defined there has a smooth fixed point, namely $(1, 0)$.
\end{pf}

\subsection*{A system of equations}

We are now ready to state and prove the main result of this section.

\begin{thm} \label{thm.oct}
Let $Q(x_1, \dots, x_m)$ be (a non-commutative and non-associative)
homogeneous polynomial of even degree in $x_1, \dots, x_m$ such that
$Q(\epsilon_1, \dots, \epsilon_m)\neq 0$ for any $(2s)$-th roots of unity
$\epsilon_1,\dots,\epsilon_m$, and let $m$ and $s$ be positive
integers. Then the system
\begin{equation} \label{e.oct-system1}
\left\{ \begin{array}{c} 
\tr(x_1^{2s}) = \dots = \tr(x_{m}^{2s}) \\ 
\tr(Q(x_1, \dots, x_m)) = 0 \, . \end{array} \right.  
\end{equation}
has no non-zero solutions in any octonion algebra of the form 
$\RMaps_{G_2}(X, \Or)$, where $X$ is a generically free $G_2$-variety
with a smooth $H$-fixed point. In particular, the system \eqref{e.oct-system1}
has no nontrivial solutions in the generic octonion algebra $\Or_{gen}$.
\end{thm}

Here $H = \lf<\tau_1, \tau_2, \tau_3 \r> \simeq (\bbZ/2\bbZ)^3$ 
is the subgroup of $G_2$ defined in~\eqref{e.oct1}.

\begin{pf}
According to Proposition~\ref{prop.oct3}, it is enough to check that
the variety
\[
Y=\Bigl\{\,(U_1:\dots:U_m)\in\bbP(\Or^m)\Bigm|
\tr(U_1^{2s}) = \dots = \tr(U_{m}^{2s}),\ \tr(Q(U_1, \dots, U_m)) = 0
\;\Bigr\}
\]
(where $U_1,\dots,U_m\in\Or$ are taken up to multiplication by an
element of $k$) has no $H$-fixed points.

A point $(U_1:\dots:U_m)\in\bbP(\Or^m)$ is $H$-fixed iff all $U_r$ lie
in the same character space for the $H$-action on $\Or$. In other
words, there exists a $\zeta \in \{ 1, i, j, l, ij, il, jl, ijl \}$
such that every $U_r$ is of the form 
$U_r=u_r\zeta$ for some $u_r \in k$. 
Note that at least one $u_r$ is non-zero; otherwise
the point $(U_1:\dots:U_m)$ is not well-defined in $\bbP(\Or^m)$. 
The condition that such a fixed point lies in $Y$ translates into 
the system 
\[
\left\{
\begin{aligned}
& u_1^{2s}=\dots=u_m^{2s}\\
& Q(u_1, \dots, u_m)=0
\end{aligned}
\right.
\]
of homogeneous equations in $u_1,\dots,u_m$. If $u_1 = 0$ then
the remaining $u_r$ are also equal to $0$, a contradiction. If $u_1 \neq 0$
then $\epsilon_r=u_r/u_1$ is a $(2s)$-th root of unity for each $r=1,\dots,m$,
and $Q(\epsilon_1, \dots, \epsilon_m)=0$, contradicting our assumption 
on $Q$. This shows that $Y$ has no $H$-fixed points.
\end{pf}

\end{document}